\documentstyle[12pt]{amsart}

\newtheorem{theorem}{Theorem}[section]

\newtheorem{definition}{Definition}[section]

\newtheorem{lemma}{Lemma}[section]

\newenvironment{proof}{\noindent {\bf {Proof:}}}{$\Box$}

\newcommand{\ZHS}{\Bbb Z{\text{-HS}}}

\numberwithin{equation}{section}

\begin{document} 

\title[Finite type invariants: A survey]{Finite type invariants of 
integral homology 3-spheres: 
A survey}
\author{Xiao-Song Lin}
\address{Department of Mathematics, University of California, Riverside,
CA 92521}
\email{xl@@math.ucr.edu}
\thanks{Research supported in part by NSF. To appear in the Proceedings 
of Warsaw Worshop in Knot Theory, July-August 1995 }
\maketitle


\section{Introduction}

We are now embarrassingly rich in knot and 3-manifold invariants. We have to 
organize these invariants systematically and find out ways to make use of 
them. The theory of finite type knot invariants, or Vassiliev invariants,
has been very successful in accomplishing the first task. Recently, 
an analogous theory of finite type invariants of integral homology 3-spheres 
started to emerge. The analogy is mainly based on the common goal of bringing 
some order to the multitude of invariants by finding some universal properties
they obey. 

If we think of quantum 3-manifold invariants of Reshetikhin and Turave 
\cite{RT} 
as non-perturbative ones, their perturbative version \cite{O1} 
is the other source
of motivation for this developing theory of finite type invariants of
integral homology 3-spheres. Non-perturbative invariants are somehow packed 
together tightly so that they usually support some 
very rich algebraic structures. 
Perturbative invariants, on the other hand, seem to be quite independent
with each other. One may see this from Theorem 3.1, the only new result in 
this paper, which claims that the space of finite type invariants of 
integral homology 3-spheres is a polynomial algebra. This leads to the 
speculation that perturbative invariants may contain more geometrical or
topological information than non-perturbative ones, at least to compensate 
the loss of algebraic richness. The recent study of Ohtsuki's perturbative
invariants in \cite{LW2} supports this speculation (see Theorem 8.1). 
    
As the title says, this is a survey of the developing theory of finite type
invariants of integral homology 3-spheres. Of course, the way we present the
material and the questions addressed depend on personal taste. But we do hope 
that we have sketched a more or less complete picture of the current status
of this theory.

This paper is based on the lectures given
in the workshop on knot theory at Banach International Center of Mathematics,
Warsaw, July 1995. We would like to thank the hospitality of Banach Center
and local organizers. We also wish to thank Zhenghan Wang for continuing
stimulating conversations, and to Stavros Garoufalidis and Jerry Levine 
for some comments.  

\section{Definitions and basic properties}

Let $\lambda$ be an invariant of oriented $\Bbb Z$-homology 3-spheres 
($\ZHS$'s in short).
We may define a difference operation on $\lambda$ with respect to a knot $K$ 
in a $\ZHS$ $M$:
$$(D^{\pm}\lambda)(M,K)=\lambda(M)-\lambda(M_{K,\pm1})$$
where $M_{K,\pm1}$ is the $\ZHS$ obtained from $M$ by a $\pm1$-surgery on $K$.
Roughly speaking, the invariant $\lambda$ is said to be of {\it finite type} 
if it vanishes under a finite iteration of the difference operation $D^{\pm}$.
See Definition 2.2 below. To be more specific, we proceed in a slightly 
different direction first.
 
An {\it algebraically split link} (ASL in short) in a 
$\ZHS$
 $M$ is a link with unoriented and unordered 
components such that all its linking numbers are zero. 
Henceforth, all 3-manifolds will be oriented $\Bbb Z$-homology 3-spheres
and all links will be $\pm1$-framed ASL's. We will denote by $\#L$ the number 
of components of a link $L$. If $L$ is a link in a $\Bbb Z$-HS $M$, 
we denote by
$M_L$ the $\Bbb Z$-HS obtained from surgery on $L$. 
In particular, if $L$ is an empty link, then $M_L=M$.

Let $\cal I$ be the set of orientation preserving 
homeomorphism classes of $\ZHS$'s. 
Let $\cal  M$ be the vector space over $\Bbb Q$ with $\cal I$ as 
a basis. For $M\in{\cal I}$ and a link $L$ in $M$, define an element $(M,L)\in
{\cal M}$ by
$$(M,L)=\sum_{L'\subset L}(-1)^{\#L'}M_{L'}$$
where $L'$ runs through all sublinks (including the empty link) of $L$. Let 
${\cal M}_k$ be the subspace of $\cal M$ spanned by $(M,L)$ for 
all $M$ and all $L$ with $\#L=k+1$. The space has a natural stratification:
$${\cal M}={\cal M}_{-1}\supset{\cal M}_0\supset\cdots\supset{\cal M}_{k-1}
\supset{\cal M}_k\supset\cdots.$$

For $A,B\in{\cal M}_{k-1}$, we denote $A\sim B$ if $A-B\in{\cal M}_k$. 
We have the following basic lemmas from \cite{O2}.

\begin{lemma} For any $(M,L)\in{\cal M}_{k-1}$, there is a link $J$ in 
$S^3$ with $\#J=k$ such that $(M,L)\sim(S^3,J)$.
\end{lemma}

\begin{lemma} $(S^3,L)\sim(S^3,J)$ if $L$ and $J$ are surgery equivalent.
\end{lemma}

\begin{lemma} Let $L$ and $J$ be the same link in $S^3$ with 
different $\pm1$ framings. Let $s(L)$ be the product of framings of $L$. 
Then $s(L)(S^3,L)\sim s(J)(S^3,J)$.
\end{lemma}

We have to recall here the definition of surgery equivalence of links from 
\cite{Le}.

\begin{definition} A surgery modification on a link $L$ in $S^3$ is first to 
have a disk $B$ in $S^3$ intersecting $L$ only in its interior with zero 
algebraic intersection number. Then perform a $\pm1$ surgery on $\partial B$ 
which will change $L$ to another link $J$ in $S^3$. Two links $L$ and $J$ are 
called {\em surgery equivalent} if one of them can be changed to another by a 
finite sequence of surgery modifications.
\end{definition}
 
From the lemmas above, we may think of ${\cal M}_{k-1}/{\cal M}_k$ as 
spanned by 
surgery equivalence classes of ASL's in $S^3$ with $k$ components, where an 
ASL $L$ in $S^3$ with $k$ components is regarded as an element $(S^3,L)\in
{\cal M}_{k-1}$ with the understanding that the framings of $L$ are all $+1$.

In \cite{Le}, Levine showed that the surgery equivalence classes of ASL's with 
oriented and ordered components are classified by Milnor's triple $\mu$-
invariants. To be more specific, surgery equivalence classes of ASL's with $k$
oriented and ordered components are in one-one correspondence with collections
of integers $\{\mu(i_1i_2i_3)\}$, where $i_1,i_2,i_3$ are distinct indices 
among $1,\dots,k$, such that
$$\left\{
\begin{aligned}&\mu(i_1i_2i_3)=\mu(i_2i_3i_1)=\mu(i_3i_1i_2)\\
&\mu(i_1i_2i_3)=-\mu(i_2i_1i_3)
\end{aligned}\right.$$
Such a collection of integers $\{\mu(i_1i_2i_3)\}$ is realized by an ASL with
$k$ oriented and order components in the following way: First, for each 
unordered triple of indices $\{i_1,i_2,i_3\}$ with $\mu(i_1i_2i_3)\neq0$,
we construct an ASL with 
three oriented components named by $i_1,i_2,i_3$, respectively, such that 
its triple Milnor $\mu$-invariant is the given one $\mu(i_1i_2i_3)$. This can 
be done by using (some variations of) the Borromean rings. Then, by 
some appropriate band sums, we connect all circles with the same name 
together to get the desired ASL. 

In \cite{O2}, Ohtsuki showed that $\text{dim}\,
({\cal M}_{k-1}/{\cal M}_k)<\infty$. The proof of this finiteness theorem  
from the point of view of Levine's classification of surgery equivalence 
classes of ASL's were given in \cite{GL1} and \cite{GrLi}. We
quote here the version in \cite{GrLi}. 

\begin{theorem} In  ${\cal M}_{k-1}/{\cal M}_k$, each ASL can be expressed as a
linear combination of ASL's whose collections of triple $\mu$-invariants
satisfy the following conditions:
\begin{enumerate}
\item each $\mu(i_1i_2i_3)$ is either $0$ or $1$ for certain fixed
cyclic order of indices;
\item each index appears in at most two non-zero $\mu(i_1i_2i_3)$'s.
\end{enumerate}
\end{theorem}

Thus, we see that ${\cal M}_{k-1}/{\cal M}_k$ has a spanning set in one-one 
correspondence with graphs having $k$ edges such that
\begin{itemize}
\item each vertex is of valence either 1 or 3;
\item each edge connects two distinct vertices;
\item no three edges have the same set of end points;
\item no edges are isolated.
\end{itemize}    

Furthermore, it was proved in \cite{GL1} that ${\cal M}_{k-1}/{\cal M}_k$ 
has a 
spanning set in one-one correspondence with graphs having $k$ edges such that
\begin{itemize}
\item each connected component is either a graph whose vertices are all of 
valence 3 or a graph with three edges and one valence 3 vertex, like the letter
Y;
\item each edge connects two distinct vertices;
\item no three edges have the same set of end points.
\end{itemize}
Thus, we have the following theorem.
 
\begin{theorem}{\em (Garoufalidis-Levine)} ${\cal M}_{3n}={\cal M}_{3n+1}=
{\cal M}_{3n+2}$.
\end{theorem}

Now let us consider rational $\Bbb Z$-HS invariants. 
Such an invariant extends naturally to a linear functional on $\cal M$.

\begin{definition}{\em (Ohtsuki \cite{O2})} A rational $\Bbb Z$-HS 
invariant $\lambda$ is said to be a {\em finite type invariant} 
of order $\leq k$ if $\lambda|{\cal M}_k=0$.
\end{definition}

Obviously, the set ${\cal O}_k$ of all finite type $\ZHS$ invariants 
of order $\leq k$ form a finite dimensional vector space over 
$\Bbb Q$. In fact, 
${\cal O}_k$ is the dual of ${\cal M}/{\cal M}_k$ by definition. As a 
corollary of Theorem 2.2, we have the following theorem.

\begin{theorem}{\em (Garoufalidis-Levine)} ${\cal O}_{3n}={\cal O}_{3n+1}=
{\cal O}_{3n+2}$.
\end{theorem}

In \cite{Ga}, a variation of Definition 2.2 was introduced. It is based on the 
notion of boundary links. A {\it boundary link} in a 3-manifold $M$ is a
link whose components bound disjoint Seifert surfaces. Obviously, a boundary
link is an ASL but not vice versa. Let ${\cal W}_n$ be the subspace of 
${\cal M}$
spanned by pairs $(M,L)$ where $L$ is a boundary link in $M$ with $n+1$
components. A $\ZHS$ invariant $\lambda$ is said 
of order $\leq n$ in the sense of Garoufalidis if $\lambda|{\cal W}_n=0$.
We denote by ${\cal G}_n$ the vector space of invariants of order $\leq n$
in the sense of Garoufalidis.

\begin{theorem}{\em (Garoufalidis-Levine \cite{GL2})} ${\cal O}_{3n}\subseteq
{\cal G}_{n}$.
\end{theorem}     

It is conjectured in \cite{GL2} that ${\cal G}_n={\cal O}_{3n}$.

Apparently, the space ${\cal G}_n$ is harder to work with than 
${\cal O}_{3n}$ and 
much less is known about it than about ${\cal O}_{3n}$. But its 
relationship with Heegaard splittings of 3-manifolds might provide some 
additional useful information about finite type invariants and the topology
of 3-manifolds.

\section{The space of finite type invariants as a polynomial algebra}

As observed in \cite{Ga}, if $\lambda\in{\cal O}_k$ and 
$\lambda'\in{\cal O}_l$, then
$\lambda\lambda'\in{\cal O}_{k+l}$. Thus, the space of all finite type $\ZHS$
invariants is a graded commutative algebra. The purpose 
of this section is to show that
the graded algebra of all finite type $\ZHS$
invariants is actually a graded commutative
and cocommutative Hopf algebra. By the structure theorem of graded commutative
and cocommutative Hopf algebras \cite{MM}, we conclude that the graded algebra
 of all 
finite type $\ZHS$ invariants is a polynomial 
algebra generated by primitive invariants. Furthermore, our construction 
of the coproduct implies that primitive invariants are exactly those 
which are additive under the connected sum of 3-manifolds. 

It might be of some interest to comment on the history of the proof, which 
we are going to present, of the fact 
that the graded algebra of finite type $\ZHS$ invariants 
is a graded commutative and cocommutative Hopf algebra. 
The same argument applied to finite
type knot invariants was first presented by this author in an AMS regional 
conference and the West Coast Topology Symposium at Stanford University 
in the spring of 1992. Before that, it was known to Bar-Natan 
\cite{BN1} 
that the graded algebra of weight systems on chord diagrams is a graded 
commutative and commutative Hopf algebra. Kontsevich's work \cite{Ko}, which 
also appeared 
in the spring of 1992, surpassed our direct argument by showing that the 
graded algebra of finite type knot invariants is isomorphic to the graded
algebra of weight systems on chord diagram (but see the remark at the end 
of this section). Therefore, our direct argument 
that the graded algebra of finite type knot invariants is a graded commutative
and cocommutative Hopf algebra never appeared in publication.   
 
In \cite{GO}, the notion of weight systems for finite type invariants
of $\Bbb Z$-homology 
3-spheres was introduced and it was shown that these weight systems form 
a graded commutative and cocommutative Hopf algebra. Since we don't know 
yet whether we could integrating every weight system in the sense of 
Garoufalidis and Ohtsuki to a $\ZHS$ invariant as in the case
of knot invariants, it is worthwhile now to argue directly that the graded
algebra of finite type $\ZHS$ invariants is also
a graded commutative and cocommutative Hopf algebra. It is obvious from the 
construction that the map from finite type invariants to their weight
systems is a graded Hopf algebra homomorphism .

To begin with, we denote by
$${\cal O}=\bigcup_{k=0}^{\infty}{\cal O}_k$$
the space of all finite type $\ZHS$ invariants.

\begin{lemma} For $\lambda\in{\cal O}_k$ and $\lambda'\in{\cal O}_l$, let
$\lambda\lambda'$ be the multiplication of $\lambda$ and $\lambda'$ as linear
functionals on ${\cal M}$, we have $\lambda\lambda'\in{\cal O}_{k+l}$.
\end{lemma}

\begin{proof} If we think of $\lambda|{\cal M}_{n-1}$ as the $n$-th difference
of $\lambda$ as in the case of knot invariants \cite{BN2},
this lemma is simply 
a corollary of the Leibniz formula for 
differences of the product of two functionals.
\end{proof}

Let $\lambda$ be a $\ZHS$ invariant. We denote by $M\#N=N\#M$ the 
connected sum of two 
oriented 3-manifolds, which inherits an orientation from that of $M$ and $N$.
Fix the factor $N$ in the connected sum and we will get a new $\ZHS$ invariant
$\lambda_N$:
$$\lambda_N(M)=\lambda(M\#N).$$

\begin{lemma} If $\lambda\in{\cal O}_k$, then $\lambda_N\in{\cal O}_k$.
\end{lemma}

\begin{proof} This comes directly from the definition of ${\cal O}_k$.
\end{proof}

We will chose a basis for $\cal O$ in the following way. Since
$${\cal O}_0\subset{\cal O}_1\subset\cdots\subset{\cal O}_{k-1}\subset
{\cal O}_{k}
\subset\cdots$$
and each ${\cal O}_k/{\cal O}_{k-1}$ is finite dimensional, there is a basis
$\{\omega_0,\omega_1,\dots,\omega_n,\dots\}$ of $\cal O$ 
such that for a certain 
sequence of non-negative
integers $n_0\leq n_1\leq\cdots\leq k_n\leq\cdots$, $\{
\omega_0,\omega_1,\dots,\omega_{n_k}\}$ is a basis for ${\cal O}_k$. 
Furthermore,
we assume that $\omega_0=1$, i.e. the constant invariant assigning to every 
$\ZHS$ the number 1, and $\omega_n(S^3)=0$ if $n>0$.
 
Now let $\lambda\in{\cal O}_k$, then $\lambda_N\in{\cal O}_k$. We write
$$\lambda_N=\phi_0(N)\omega_0+\phi_1(N)\omega_1+\cdots+\phi_{n_k}(N)\omega_
{n_k}$$
where the coefficients $\phi_0(N)$, $\phi_1(N)$, $\dots$, and $\phi_{n_k}(N)$
can be thought of as $\ZHS$ invariants. 

\begin{lemma} The $\ZHS$ invariants $\phi_0,\phi_1,\dots,\phi_{n_k}$ are of 
finite type. More precisely, if the order of $\omega_l$ is $\leq r$, then the 
order of $\phi_l$ is $\leq k-r$.
\end{lemma}

\begin{proof} The first conclusion is very easy to obtain: Simply notice that
$\omega_0,\omega_1,\dots,\omega_{n_k}$ are linearly independent and 
$\lambda_N(M)=\lambda(M\#N)$ by definition. Since $\lambda\in{\cal O}_k$,
for every ASL $L$ in $N$ with $k+1$ components, we have
$$\lambda_{(N,L)}(M)=\lambda(M\#(N,L))=0$$
for all $M$. This implies that 
$$\phi_0((N,L))=\phi_1((N,L))=\cdots=\phi_{n_k}((N,L))=0.$$

The second conclusion is obtained using the same argument, but taking into
consideration the way we chose the basis $\{\omega_0,\omega_1,\dots,
\omega_{n_k}\}$. We will consider links $J$ in $M$ and $L$ in $N$ such
that $\#J+\#L=k+1$. We have
$$\lambda_{(N,L)}((M,J))=\lambda((M\#N,J\coprod L))=0.$$
Then, we may argue as before to finish the proof.
\end{proof}

We define a coproduct $\Delta:{\cal O}\rightarrow{\cal O}\otimes{\cal O}$ 
preserving
the grading by
$$\Delta(\lambda)=\omega_0\otimes\phi_0+\omega_1\otimes\phi_1+\cdots+
\omega_{n_k}\otimes\phi_{n_k}.$$

\begin{lemma} The coproduct defined above is cocommutative. Moreover, if 
order$\,(\lambda)>0$, we have
$$\Delta(\lambda)=1\otimes\lambda+\lambda\otimes1+\sum_{i,j}\alpha_i
\otimes\beta_j$$
with homogeneous elements $\alpha_i,\beta_j$ of positive orders such that
$\text{order}\,(\alpha_i)+\text{order}\,(\beta_j)=\text{order}\,(\lambda)$.
\end{lemma}

\begin{proof} The first conclusion is because of the commutativity 
of the connected sum of 
3-manifolds.

For the second conclusion, we first get $\phi_0=\lambda$ by taking $M=
S^3$. So, we get the summand $1\otimes\lambda$ in $\Delta(\lambda)$. The 
summand $\lambda\otimes1$ then comes from the cocommutativity of $\Delta$.
\end{proof}

A $\ZHS$ invariant $\lambda$ is called {\it additive} if
$$\lambda(M\#N)=\lambda(M)+\lambda(N).$$ 
A $\ZHS$ invariant $\lambda\in{\cal O}$ is called {\it primitive} if 
$$\Delta(\lambda)=1\otimes\lambda+\lambda\otimes1.$$

\begin{lemma} $\lambda\in{\cal O}$ is additive iff it is primitive.
\end{lemma}

\begin{proof} This directly follows from the definition of the coproduct
$\Delta$.
\end{proof}

Now we may summarize the previous discussion into the following theorem.

\begin{theorem} $\cal O$ is a polynomial algebra generated by additive 
invariants.
\end{theorem}

\begin{proof} By the previous lemmas, $\cal O$ is a graded commutative 
and cocommutative Hopf algebra. By the structure theorem of graded 
commutative and cocommutative Hopf algebras \cite{MM}, we know $\cal O$ is a 
polynomial algebra generated by primitive invariants. The theorem follows 
since primitive invariants are exactly additive invariants.
\end{proof} 

\medskip

\noindent{\bf Remark:} Instead of defining $\cal M$ as a vector space over
a field (say $\Bbb Q$) spanned by all $\Bbb Z$-HS's, we may define it as the
free $R$-module generated by all $\Bbb Z$-HS's, where $R$ is a commutative ring
with unit 1. The subspaces ${\cal M}_k$ then become submodules. 
As pointed out by J\`ozef Przytycki, similar to the case of knots discussed in 
\cite{Pr}, the completion of $\cal M$ with respect to the stratification 
$\{{\cal M}_k\}$ is a commutative and cocommutative $R$-Hopf algebra. 
Again similar to the case of knots, the 
obvious question is 
whether the modules ${\cal M}_k/{\cal M}_{k+1}$ are
torsion free. 

\section{Basic examples: Ohtsuki's invariants}
 
We are mainly talking about generalities so far and haven't had any non-trivial
finite type $\ZHS$ invariants yet. Our basic examples come from Ohtsuki's work 
\cite{O1}. In \cite{O1}, Ohtsuki extracted a series of 
rational $\ZHS$ invariants from 
the $SU(2)$ quantum invariants of Reshetikhin and Turaev \cite{RT}. 
Physically, they 
correspond to the coefficients of the asymptotic expansion of Witten's 
Chern-Simons path integral at the trivial connection as shown by Rozansky
\cite{R1,R2}.
Rozansky has also argued, based on some physical considerations, that Ohtsuki's
invariants are all of finite type.

Let first us explain briefly how Ohtsuki's invariants are derived. 

Let $\tau_r(M)$ be the $SU(2)$ quantum invariant of Reshetikhin and Turaev
at the $r$-th root of the unit $q=e^{2\pi\sqrt{-1}/r}$, as normalized
in \cite{KMe}. From \cite{Mu}, we know that if $M$ is a $\ZHS$,
$\tau_r(M)\in{\Bbb Z}[q]$ when $r$ is an odd prime.
So we may write
$$\tau_r(M)=a_{r,0}+a_{r,1}(q-1)+\cdots+a_{r,n}(q-1)^n+\cdots$$
for $a_{r,n}\in{\Bbb Z}$. We will fix $n$ and think of $a_{r,n}$ as a function
of $r$. We ask whether there is a rational number $\lambda_n$, 
independent of $r$, such that
\begin{equation}
a_{r,n}\equiv\lambda_n\qquad\text{mod}\,\, r
\end{equation}
is true for all odd primes $r$ sufficiently large. 
Of course, since $q$ is not an 
indeterminant, $a_{r,n}$ here is not well defined. Nevertheless, as $q^r=1$,
$a_{r,n}$ is well defined modulo $r$ and our question is thus well posed.

The situation here is rather like the one dealt with in the so-called Fermat's
little theorem: Fix an integer $a$, $a^{r-1}\equiv 1$ mod $r$ 
for all odd primes
$r$ sufficiently large. From elementary number theory \cite{HW}, 
we know that things
are not always so nice. For example, $(\frac{r-1}{2})!$ mod $r$ depends on $r$.
For this reason, we call in \cite{LW2} the formal power series
$$\sum_{n=0}^\infty\lambda_n(t-1)^n$$
the {\it Fermat limit} of $\tau_r(M)$ if (4.1) is true for every $n$, 
and we denote
$$\text{f-lim}\,\tau_r(M)=\sum_{n=0}^\infty\lambda_n(t-1)^n.$$
Certainly, if the Fermat limit of $\tau_r(M)$ exists, it is unique.

\begin{theorem}{\em (Ohtsuki \cite{O1}, see also \cite{LW2})} 
The Fermat limit of $\tau_r(M)$
exists. In particular, we get a sequence of rational $\ZHS$ invariants
$\lambda_n=\lambda_n(M)$.
\end{theorem}

It is quite easy to get $\lambda_0=1$. If we denote by $\lambda_{C}$ the Casson
invariant of $\ZHS$'s \cite{AM}, we have the following theorem.

\begin{theorem}{\em (Murakami \cite{Mu})} $\lambda_1=6\lambda_C$.
\end{theorem}

The first example of non-trivial finite type $\ZHS$ invariants 
comes from $\lambda_C$, or equivalently, $\lambda_1$. 
We combine many known facts into the following 
theorem. See \cite{O2,Ga}.

\begin{theorem} We have
\begin{enumerate}
\item ${\cal O}_0={\cal O}_1={\cal O}_2\cong{\Bbb Q}$  
and it is spanned by the 
constant invariant 1;
\item ${\cal O}_3={\cal O}_4={\cal O}_5\cong{\Bbb Q}^2$ and 
it is spanned by 1 and 
$\lambda_C$.
\end{enumerate}
\end{theorem}

Following Theorem 2.3, we only need to explain why $\lambda_C$ is of order 
$\leq 3$. We proceed using the surgery formula of 
$\lambda_C$ given in \cite{Ho}.

Let $\nabla(L;z)$ be the Conway polynomial \cite{Ka} 
of a link $L$ in $S^3$. Then,
$$\nabla(L;z)=z^{\#L-1}(1+a_2(L)z^2+a_4(L)z^4+\cdots+a_{2n}(L)z^{2n}).$$

If $L$ be a $\pm1$-framed ASL, we will denote by $f_L$ the product of framings
of $L$. 

\begin{theorem}{\em (Hoste \cite{Ho})} Let $L$ be a 
$\pm1$-framed ASL in $S^3$, 
then 
$$\lambda_C(S^3_L)=\sum_{L'\subset L}f_{L'}a_2(L').$$
\end{theorem}
 
The following is a very interesting property of the link invariant
$a_2$ on ASL's.

\begin{lemma} {\em (Hoste \cite{Ho})} If $L$ is an ASL with $\#L>3$, then
$a_2(L)=0$.
\end{lemma}

We will leave it as an exercise for the reader to show that 
$\lambda_C$ is of order $\leq 3$ using Theorem 4.4 and Lemma 4.1. 
Since $\lambda_C$ is
not a constant invariant, it is of order 3. 

It seems to be a common belief that $\lambda_n$ should be of order
$3n$ in general, and this has been verified physically by Rozansky 
\cite{R1,R2}.  
Surgery formulae for $\lambda_n$ were given in \cite{LW2}. Unfortunately, some 
coefficients in these formulae can only be determined recursively. 
We don't have close formulae yet to express these coefficients appeared
in the surgery formulae for $\lambda_n$ when $n>2$. This makes it difficult to
prove that $\lambda_n$ is of order $3n$ for $n>3$. On the other hand, the 
explicit surgery formula for $\lambda_2$ in \cite{LW2} 
makes it possible to prove
that $\lambda_2$ is of order 6.

\begin{theorem}{\em (Lin-Wang \cite{LW3})} $\lambda_2$ is of order 6. 
\end{theorem}

We will briefly explain the proof of this theorem in the next section when
the surgery formula for $\lambda_2$ is given.

In \cite{GO}, it is estimated that $\text{dim}\,{\cal O}_6\leq 4$ by using 
some relations among trivalent graphs representing $\ZHS$'s. 
On the other hand, a simple computation using the 
surgery formula for $\lambda_2$ (see (5.1) below) shows that
$\lambda_1^2$ and $\lambda_2$ are linearly independent order 6 
invariants. Thus,
$${\cal O}_6={\cal O}_7={\cal O}_8\cong{\Bbb Q}^4$$
and it is spanned by 
$1,\lambda_1,\lambda_1^2, \lambda_2$.

\section{The surgery formula for $\lambda_2$}

We will denote by $O$ the unknot and $\emptyset$ the empty link. 

The Jones polynomial $V(L;t)\in{\Bbb Z}[t^{\pm\frac12}]$ of a link $L$ in $S^3$
is defined by
$$\left\{\begin{aligned}
&V(O;t)=1;\\
&V(\emptyset;t)=(t^{\frac12}+t^{-\frac12})^{-1};\\
&tV(L_+;t)-t^{-1}(L_-;t)=(t^{\frac12}-t^{-\frac12})V(L_0;t).
\end{aligned}
\right.
$$ 
Here, as usual, $L_+$, $L_-$ and $L_0$ are three links which have plane 
projections identical to each other except in one small disk where their 
projections are a positive crossing, a negative crossing and an orientation
smoothing of that crossing, respectively. Note that the Jones polynomial we 
use here differs from the usual definition in \cite{J} 
with a change of variable
$t$ to $t^{-1}$ and a factor of $(-1)^{\#L-1}$. 

We put
$$X(L;t)=\frac{V(L;t)}{(t^{\frac12}+t^{-\frac12})^{\#L-1}}$$
and
$$\Phi(L,t)=\sum_{L'\subset L}(-1)^{\#L-\#L'}X(L';t)$$
where the sum runs through all sublinks of $L$. We need to use the derivatives 
of $\Phi(L;t)$ at $t=1$:
$$\Phi_i(L)=\left.\frac{d^i\Phi(L;t)}{dt^i}\right|_{t=1}.$$
The basic link invariants
used in formulae for $\lambda_n$'s are
$$\phi_i(L)=\frac{(-2)^{\#L}}{(\#L+i)!}\Phi_{\#L+i}(L)$$
for $i=1,2,\dots$ (see \cite{LW2}). According to Murakami \cite{Mu}, we have
$$\phi_1(L)=6a_2(L)$$
if $L$ is an ASL.

Given a link $L$ and a positive integer $m$, we denote by $L^m$ the 0-framed
$m$-parallel of $L$, i.e. each component in $L$ is replaced by $m$ parallel
copies having linking number zero with each other. So $L^m$ is an ASL if 
$L$ is. Sublinks of $L^m$ will be assumed to be in one-one correspondence 
with $\#L$-tuples $(i_1,i_2,\dots,i_{\#L})$, in such a way that the
corresponding sublink $L'$ will have $i_\xi$ parallel copies of the $\xi$-th 
component of $L$, $0\leq i_\xi\leq m$. If $L$ is a framed link, $L^m$ and all 
its sublinks will inherit a framing from that of $L$.

\begin{theorem} {\em (Lin-Wang \cite{LW2})} 
Let $L$ be a $\pm1$-framed ASL. Then
\begin{equation}
\lambda_2(S^3_L)=\sum_{L'\subset L}\phi_1(L')f_{L'}\frac{\#L'}2
+\sum_{L'\subset L^2}\phi_2(L')f_{L'}\frac1{2^{s_2(L')}}.
\end{equation}
Here, if $L'\subset L$ corresponds to the $\#L$-tuple $(i_1,i_2,\dots,
i_{\#L})$, $s_2(L')=\#\{i_\xi\,;\,i_\xi=2\}$.
\end{theorem} 

Similar to the proof that $\lambda_C$ is of order $\leq3$, the following lemma
together with Theorem 5.1 implies that $\lambda_2$ is of order $\leq6$. A 
little more computation will then establish Theorem 4.5.

\begin{lemma} {\em (Lin-Wang \cite{LW3})} If $L$ is an ASL with $\#L>6$,
then $\phi_2(L)=0$
\end{lemma}

The proof of this lemma is based on studying the colored Jones polynomial 
\cite{RT}. A generalization of (the first part of) the so-called 
Melvin-Morton conjecture \cite{MeMo,BNG} for ASL's turns out to be crucial.

\section{Induced knot invariants}

Let $\lambda$ be a $\ZHS$ invariant. It induces a knot invariant in the
following way: For a knot $K$ in $S^3$, we assign to it the framing 1. 
Then we define a knot invariant $\psi_\lambda$ by
$$\psi_\lambda(K)=\lambda(S^3_K).$$

\begin{lemma} If $\lambda$ is of finite type, so is $\psi_\lambda$ as a knot 
invariant.
\end{lemma}

This is quite obvious: A crossing change crossings can be accomplished
by a $\pm1$-surgery on an appropriately positioned unknot. In fact,
the proof of this lemma shows that if $\lambda$ is of order $\leq k$, then
$\psi_\lambda$ is of order $\leq k-1$. The fact that the order of a finite
type $\ZHS$ invariant must be a multiple of 3 (Theorem 2.3) somehow implies 
a stronger estimate on the order of the induced knot invariant. 

\begin{theorem} {\em (Habegger \cite{Ha}, see also \cite{GL2})} 
If $\lambda$ is a $\ZHS$ invariant
of order $3n$, $\psi_\lambda$ is a knot invariant of order $\leq 2n$. 
\end{theorem}

For example, the induced knot invariant of $\lambda_C$ is $a_2$, i.e. the
second coefficient of the Conway polynomial. The induced knot invariant
of $\lambda_2$ is slightly more complicated but can be obtained using
(5.1). To express it using classical knot invariants, we denote
$$v_i(L)=\left.\frac{d^iV(L;e^h)}{dh^i}\right|_{h=0}.$$
Let $\psi_2$ be the knot invariant induced by $\lambda_2$, we have
\begin{equation}\psi_2(K)=\frac13\,v_2(K)-\frac13\,v_3(K)-\frac16\,v_4(K)+
\frac23\,v_2^2(K).
\end{equation}

\section{More on Ohtsuki's invariants: their integrality}

{\it A priori}, Ohtsuki's invariants $\lambda_n$ are rational invariants. From
Theorem 4.2, we see that $\lambda_1\in 6{\Bbb Z}$. 
The following theorem follows
from the surgery formula (5.1) for $\lambda_2$ and a careful study of
the coefficients of the colored Jones polynomial of ASL's.

\begin{theorem} {\em (Lin-Wang \cite{LW2})} $\lambda_2\in3{\Bbb Z}$.
\end{theorem}

In \cite{LW2}, 
it is conjectured that $n!\,\lambda_n\in 6{\Bbb Z}$. This conjecture 
is certainly motivated by Theorems 4.2 and 7.1.
But we do have some other evidence supporting this conjecture.
For example, it is proved in \cite{LW2} 
that the denominator of $\lambda_n$ can 
always be factored by some powers of $2,3,\dots,n$.
 
\section{Some applications of Ohtsuki's invariants}

Formula (5.1) allows us to compute $\lambda_2$ quite easily. We may then 
use the invariant $\lambda_2$ to distinguish various $\ZHS$'s. For example,
if $M$ is the Poincar\'e homology 3-sphere $\Sigma(2,3,5)$ ($+1$-surgery on 
the 
right-handed trefoil knot), we have $\lambda_2(M)=39$. And if $M$ is the 
homology 3-sphere $\Sigma(2,3,7)$ ($+1$-surgery on the left-handed
trefoil knot), we have $\lambda_2(M)=63$. In general, we may 
use $\lambda_2$ to distinguish $\ZHS$'s
 obtained from $1/n$-surgeries on a knot $K$ and its mirror image 
$K^*$, respectively, provided that $v_3(K)\neq0$. 

\begin{theorem} {\em (Lin-Wang \cite{LW2})} If $v_3(K)\neq0$, then
$S^3_{K,1/n}\neq S^3_{K^*,1/n}$ for each $n\neq0$.
\end{theorem}

Notice that these two $\ZHS$'s $S^3_{K,1/n}$ and $S^3_{K^*,1/n}$ can not
be distinguished using the Casson invariant.

Needless to say, we don't know what kind of geometrical or topological 
obstruction the invariant $\lambda_2$ represents which prevents these two
$\ZHS$'s from being homeomorphic.

See \cite{LW2} for more on this kind of applications of Ohtsuki's invariants. 

\section{Open questions}

Compare with the well developed theory of finite type knot invariants \cite
{BN2,BL,Gu,Li,Ko,St1,V}, much less have been done in the study of finite type 
$\ZHS$ invariants. Many questions have been asked in the literature of
this subject. We will address a few questions here which are either not asked
before or, in the case they are, picked up again because of their 
significance in our point of view.

\noindent{\bf Question 1.} Is there any other justification for Definition
2.2 besides its analogy to the definition of finite type knot invariants
and the fact that it works for the coefficients of the asymptotic expansion
of Witten's Chern-Simons path integral at the trivial connection \cite{R1,R2}?

One of the original sources of the definition of finite type knot invariants
is the study of the topology of the singularities in some functional spaces
\cite{V}. See also \cite{St1} and \cite{Li}
for treatments using only elementary topology in
3-dimension. Such a topological interpretation provides us the most natural 
explanation of the 4-term relation and why it alone suffices essentially
to define a \lq\lq weight system''. We don't have such a picture yet for
the definition of finite type $\ZHS$ invariants.

\noindent{\bf Question 2.} Does ${\cal M}/{\cal  M}_k$ constitute an abelian 
group under the connected sum? What is the most natural operation which
changes a $\ZHS$ to another one preserving the value 
of any finite type $\ZHS$ invariant up to a fixed order? If we have such an
operation on $\ZHS$'s, is there any topological property of $\ZHS$'s which 
is preserved by this operation?

The first question here is certainly motivated by Gusarov's work \cite{Gu} 
and the 
second by, among others, say, Stanford's work \cite{St2}. The operation on 
knots defined by Stanford is 
to insert a pure braid commutator somewhere in a given knot. Notice that 
Gusarov has reported (Oberwolfach, September 1995) that Stanford's operation
actually generates all knot with the same value for any finite type invariant 
up to a fixed order. Finally, the answer to the last question above 
seems to be not known even for Stanford's operation on knots. Because of
its significance, it seems to be worthwhile to repeat it again: 
Is there any topological
property of knots which is preserved by Stanford's operation?

\noindent{\bf Question 3.} Does $\lambda_2/3$ count algebraically any 
geometrical or topological objects related with the $\ZHS$ in question?

This is certainly motivated by the fact that $\lambda_1/6$ is the Casson
invariant, which counts algebraically the number of conjugacy classes
of irreducible representations into $SU(2)$ of the fundamental group
of the $\ZHS$ in question \cite{AM}. Probably, the only way to answer 
this question
is to define some counting invariants for $\ZHS$'s first and then identify it
with $\lambda_2$ using the surgery formula (5.1).

Notice that the geometric definition of the Casson invariant together
with its surgery formula leads to a criterion for detecting knots with 
property P. An answer to Question 3 will provide us another such criterion.

\noindent{\bf Question 4:} Is it true that if $\lambda$ is a $\ZHS$ invariant
of order $\leq k$, then there is a constant $C$ such that for 
every ASL $L$ in $S^3$ with $\mu$ components, we have 
$$|\lambda(S^3_L)|\leq C\mu^k?$$
Does Ohtsuki's invariants $\lambda_n$ have this property? 

A similar estimate for finite type knot invariants is conjectured in 
\cite{LW1} 
and proved in \cite{BN3}. The question to $\lambda_n$, with some additional
control over the constants $C$ involved, is probably related with
the convergence problem of the formal power series
$$\sum_{n=0}^\infty\lambda_n(t-1)^n.$$
See \cite{L} for examples of 
3-manifolds for which this formal power series does
converge.   

\section{New developments}

There have been some major new developments since the submission of this 
article in October, 1995. In what follows, for the reader's convenience, 
we give brief summaries of these new works.
 
\subsection{} In an important work of Le, Murakami and Ohtsuki \cite{LMO},
a universal  3-manifold invariant $\Omega(M)$ is 
constructed, which takes values in the 
completed, graded 
commutative and 
cocommutative Hopf algebra of trivalent graphs 
subject to the $AS$ and $IHX$
relations. To start with, let us denote by ${\cal A}(X)$ 
the graded vector space generated by graphs with univalent and trivalent 
vertices supported on $X$ where
\begin{itemize}
\item $X$ can be the empty space $\emptyset$, a set of $m$ ordered points
(which will be denoted simply by $m$), and 
the disjoint 
union of oriented circles $\coprod S^1$;
\item all univalent vertices are on $X$;
\item when $X$ is a set of $m$ 
ordered points, points in $X$ are all we have for univalent vertices, and when
$X=\coprod S^1$, there should be no separate subgraphs with only trivalent 
vertices;
\item a cyclic ordering of the edges at every trivalent vertex is given;
\item the grading is given by a half of the number of vertices;
\item linear combination of graphs are subject to the $AS$ and $IHX$ relations,
and the $STU$ relations if $X=\coprod S^1$ (see \cite{BN2}).
\end{itemize}
Denote by $\hat{\cal A}(X)$ the completion of ${\cal A}(X)$. 

With these said, the construction in \cite{LMO} proceeds as follows:

\noindent{\em Step1 .} Construct elements $T_m^n\in{\cal A}^{(m-n)}(m)$,
which are invariant under cyclic permutations of the $m$ supporting points
and characterized essentially by a certain kind of ``crossing change 
formulae''.

\noindent{\em Step 2.} For a framed link $L$, 
renormalize the universal invariant
$\hat Z(L)$ coming from Kontsevich's construction \cite{LM} to get an
invariant $\breve Z(L)\in\hat{\cal A}(\coprod S^1)$, which is a ``group-like
element''.

\noindent{\em Step 3.} Define a linear map
$$\iota_n:{\cal A}^{(k)}(\coprod S^1)\longrightarrow{\cal A}^{(k-n)}
(\emptyset)$$
by ``installing'' $T_{2k}^n$ onto each circle so that univalent vertices
on both sides match in cyclic order, and then drop off the circles. This 
extends to a linear map
$$\iota_n:\hat{\cal A}(\coprod S^1)\longrightarrow\hat{\cal A}(\emptyset).$$

\noindent{\em Step 3.} $\iota_n(\breve Z(L))\in\hat{\cal A}(\emptyset)$,
modulo $\hat{\cal A}^{(> n)}(\emptyset)$, turns out
to be independent of the orientation of $L$ and invariant
under the type II Kirby move on the framed link $L$. So, an appropriate
normalization of $\iota_n(\breve Z(L))\in\hat{\cal A}(\emptyset)$ modulo
$\hat{\cal A}^{(> n)}(\emptyset)$ gives us an invariant $\Omega_n(M)$ of the
3-manifold $M$ obtained by surgery on $L$.

\noindent{\em Step 4.} The degree $<n$ parts of $\Omega_n(M)$ turn out to be
determined by the top degree $k$ part of $\Omega_k(M)$ denoted by 
$\Omega_k^{(k)}(M)$, for all $k<n$. So, finally, let
$$\Omega(M)=1+\sum_{n=1}^\infty\Omega_n^{(n)}(M)\in\hat{\cal A}(\emptyset).$$
It turns out that $\Omega(M)$ is also a ``group-like element''. 
 
More specifically, we have $\Omega_{n+1}^{(\leq n)}(M)=d\,\Omega_n(M)$, 
where $d$ equals the cardinality of $H_1(M;{\Bbb Z})$ if the 
first Betti number of $M$ 
is 0 and $d=0$ otherwise. So, if $M$ is a $\Bbb Z$-HS,
$$\Omega(M)=1+\sum_{n=1}^\infty\Omega_n(M).$$
When both $M_1$ and $M_2$ are $\Bbb Z$-HS's, we have
$$\Omega(M_1\sharp M_2)=\Omega(M_1)\times\Omega(M_2)$$

It was proved in \cite{LMMO} that if $M$ is a $\Bbb Z$-HS, $\Omega_1(M)$ 
is essentially the Casson invariant.

We remark that this construction of Le, Murakami and Ohtsuki has the similar 
feature as the construction of $SU(2)$ quantum invariants given by
Lickorish \cite{Lo}. The role played by the elements $T_m^n$ in the 
Le-Murakami-Ohtsuki construction is very much like the role played by
the so-called Jones-Wenzl projectors in Lickorish's construction. Are they 
related in a certain way?

\subsection{} In \cite{Le'}, Le explored further the universal invariant 
$\Omega(M)$ in the case when $M$ is a $\Bbb Z$-HS and made the 
combinatorial aspect of the theory of
finite type invariants of integral homology 3-spheres completely parallel
to the theory of finite type invariants of knots. Main results of
\cite{Le'} include (all manifolds below are $\Bbb Z$-HS's):
\begin{itemize}
\item $\Omega_n(M)$ is an invariant of order $\leq 3n$. This implies, via
\cite{GO}, that the algebra of finite type invariants is isomorphic to the
algebra of weight systems, where a {\em weight system} is a linear functional
on ${\cal A}^{(n)}(\emptyset)$.
\item An operation on 
$\Bbb Z$-HS's is defined which does not alter the value any invariant
of order $\leq n$. This operation is a straightforward generalization of
Stanford's operation. 
\item ${\cal M}/{\cal M}_k$ 
constitutes an abelian group under the
connected sum. 
\end{itemize}
So parts of Question 2 in Section 9 have been answered.

\subsection{} In \cite{GO'}, the notion of finite type invariants of
rational homology 3-spheres ($\Bbb Q$-HS's) is introduced. It is shown that 
the space of
finite type $\Bbb Q$-HS invariants is a subspace of finite type $\Bbb Z$-HS
invariants. In particular, the Casson-Walker invariant \cite{Wa} of 
$\Bbb Q$-HS's is shown to be of order 3 using its surgery formula. In
\cite{Le'}, Le also proved that $\Omega_n(M)$ is of order $\leq 3n$ for 
$\Bbb Q$-HS's. 

\subsection {} The recent work of Garoufalidis and Levine \cite{GL3} provides
another perspective to the theory of finite type $\Bbb Z$-HS invariants. 
It relates the stratified vector space $\cal M$ with 
the group algebras of some subgroups of the mapping class group completed by 
powers of their argumentation ideals. The subgroups of the mapping class
group are those which are closely related with integral homology 3-spheres.
For example, the Torelli group and the subgroup of the mapping class group 
generated by Dehn surgeries on bounding 
curves on surfaces (the main result is stated differently though in 
these two cases). Recall that the structure of the latter subgroup has been 
used by S. Morita to recover the Casson invariant \cite{Mo1,Mo2}. So, 
\cite{GL3} will probably give us some hints as of how to generalize 
Morita's work 
to recover Ohtsuki's invariant $\lambda_2$ (see Question 3 of Section 9).

\end{document}